\documentclass[12pt]{article}
\usepackage{amsmath, amsthm, amsfonts}
\usepackage{amssymb}
\usepackage{amscd}
\usepackage{verbatim}

\begin{document}
\newcommand{\eqdef}{\stackrel{{\rm def}}{=}}
\newcommand{\M}{{\mathcal M}}
\newcommand{\loc}{{\mathrm{loc}}}
\newcommand{\dx}{\,\mathrm{d}x}
\newcommand{\core}{C_0^{\infty}(\Omega)}
\newcommand{\sob}{W^{1,p}(\Omega)}
\newcommand{\sobloc}{W^{1,p}_{\mathrm{loc}}(\Omega)}
\newcommand{\merhav}{{\mathcal D}^{1,2}(\R^N)}
\newcommand{\be}{\begin{equation}}
\newcommand{\ee}{\end{equation}}
\newcommand{\mysection}[1]{\section{#1}\setcounter{equation}{0}}
%%%%%%%%%%%%%%%
\newcommand{\bea}{\begin{eqnarray}}
\newcommand{\eea}{\end{eqnarray}}
\newcommand{\bean}{\begin{eqnarray*}}
\newcommand{\eean}{\end{eqnarray*}}
\newcommand{\thkl}{\rule[-.5mm]{.3mm}{3mm}}
%%%%%%%%%%%%%%%%%%%%%%%%%%%
\newcommand{\cw}{\stackrel{\rightharpoonup}{\rightharpoonup}}
\newcommand{\id}{\operatorname{id}}
\newcommand{\supp}{\operatorname{supp}}
\newcommand{\wlim}{\mbox{ w-lim }}
\newcommand{\mymu}{{x_N^{-p_*}}}
\newcommand{\R}{{\mathbb R}}
\newcommand{\N}{{\mathbb N}}
\newcommand{\Z}{{\mathbb Z}}
\newcommand{\Q}{{\mathbb Q}}
\newcommand{\abs}[1]{\lvert#1\rvert}
%%%%%%%%%%%
\newtheorem{theorem}{Theorem}[section]
\newtheorem{corollary}[theorem]{Corollary}
\newtheorem{lemma}[theorem]{Lemma}
\newtheorem{definition}[theorem]{Definition}
\newtheorem{remark}[theorem]{Remark}
\newtheorem{proposition}[theorem]{Proposition}
\newtheorem{assertion}[theorem]{Assertion}
\newtheorem{problem}[theorem]{Problem}
%%%%%%%%%%%%%%%%%%
\newtheorem{conjecture}[theorem]{Conjecture}
\newtheorem{question}[theorem]{Question}
\newtheorem{example}[theorem]{Example}
\newtheorem{Thm}[theorem]{Theorem}
\newtheorem{Lem}[theorem]{Lemma}
\newtheorem{Pro}[theorem]{Proposition}
\newtheorem{Def}[theorem]{Definition}
\newtheorem{Exa}[theorem]{Example}
\newtheorem{Exs}[theorem]{Examples}
\newtheorem{Rems}[theorem]{Remarks}
\newtheorem{Rem}[theorem]{Remark}

\newtheorem{Cor}[theorem]{Corollary}
\newtheorem{Conj}[theorem]{Conjecture}
\newtheorem{Prob}[theorem]{Problem}
\newtheorem{Ques}[theorem]{Question}
\newcommand{\pf}{\noindent \mbox{{\bf Proof}: }}

%%%%%%%%%%%%%%%%%%
%\newenvironment{proof}{{\bf Proof.}}{\hfill $\bowtie$\vskip4mm}

\renewcommand{\theequation}{\thesection.\arabic{equation}}
\catcode`@=11 \@addtoreset{equation}{section} \catcode`@=12

\title{Concentration-compactness at the mountain pass level in semilinear elliptic problems}
\author{Kyril Tintarev\thanks{Research partly done while visiting Ceremath - University of Toulouse 1.
}
\\{\small Department of Mathematics}\\{\small Uppsala University}\\
{\small SE-751 06 Uppsala, Sweden}\\{\small
kyril.tintarev@math.uu.se}}
%\date{}
\maketitle

\begin{abstract}
\begin{small}The concentration compactness framework for semilinear elliptic equations without compactness, set originally by P.-L.Lions for constrained minimization in the case of homogeneous nonlinearity, is extended here to the case of general nonlinearities in the standard mountain pass setting of Ambrosetti--Rabinowitz. In these setting, existence of solutions at the mountain pass level $c$ is verified under a single assumption $c<c_\infty$, where $c_\infty$ is the mountain pass level for the asymptotic problem, which is completely analogous to the Lions' case. Problems on $\R^N$ and problems with critical nonlinearity are considered. Particular attention is given to nonhomogeneous critical nonlinearities that oscillate about the ``critical stem" $F(u)=|u|^{2^*}$.
\end{small}
\end{abstract}

%\end{titlepage}
%%%%%%%%%%%%%%%%%%%%%%%%%%%%%%%%%%%%%%%%%%%%%%%%%

\section{Introduction}
In this paper we study existence for the classical
semilinear elliptic problem in a domain $\Omega\subset\R^N$, 
\be \label{Theeq} -\Delta
u+\lambda u= f(x,u), \ee
with the Dirichlet boundary condition. The number $\lambda$ is assumed to be greater than the bottom of the spectrum for the Dirichlet Laplacian in $\Omega$. We consider here cases where the correspondent Sobolev imbedding lacks compactness, namely when $\Omega=\R^N$ or when the nonlinearity $f(x,s)$ has the growth of the critical magnitude $|s|^{2^*-1}$, where $2^*=\frac{2N}{N-2}$ and $N>2$. The variational framework set by Ambrosetti and Rabinowitz in \cite{AR}, initially for the case of bounded $\Omega$ and the subcritical nonlinearity, where one can benefit from compactness of Sobolev imbeddings, faces significant difficulties when the the problem lacks compactness, requiring the concentration compactness argument.
 
\subsection{Problems with homogeneity}

If $f(x,u)=a(x)|u|^{p-2}u$, solutions to \eqref{Theeq}, can be obtained by constrained minimization
\be
\label{isop}
c=\inf_{u\neq 0}\dfrac{\int_\Omega(|\nabla u|^2+\lambda u^2)dx}{\left(\int_\Omega a(x)|u|^pdx\right)^{2/p}}.
\ee
Homogeneity of the nonlinearity means that one can relplace the Lagrange multiplier for the minimizer with $1$ by multiplying the minimizer of \eqref{isop} by an appropriate positive constant. Since the loss of compactness in these problems is due to non-compact transformations (shifts or dilations), discrete sequences of these transformation define 
asymptotic problems (or problems at infinity). Applying such sequences one can immediately see that the inequality $c\le c_\infty$, where $c_\infty$ is the constrained minimum for an asymptotic counterpart of \eqref{isop}, is always true. In the famous four-paper series of P.L.-Lions (\cite{PLL1a,PLL1b,PLL2a,PLL2b}, the {\em strict} inequality $c< c_\infty$ was employed as a sufficient existence condition to prove existence minimizers in \eqref{isop} both in the subcritical ($p<2^*$) and the critical ($p=2^*$, in particular the famous Brezis--Nirenberg case $\lambda<0$, $a=1$, $N>2$, $p=2^*$ on bounded domains, \cite{BN}). 

In the autonomous ($f(x,s)=f(s)$) subcritical case on $\R^N$, the imbedding of the subspace of $H^1$- radially symmetric functions into $L^p$ is compact (\cite{Esteban}), so that existence of minimizers in \eqref{isop} follows from a standard weak continuity/ lower semicontinuity argument. Lagrange multipliers here can be set to $1$ even for general nonlinearities $f(s)$ of subcritical growth, since autonomous problems on $\R^N$ possess additional homogeneity, namely the one with respect to the transformation $s\mapsto u(\cdot/s)$. Existence results for the subcritical autonomous case are due to Berestycki and Lions\cite{BerestickyLions}. In the critical case the imbedding of the subspace of radial functions into $L^{2^*}$ is no longer compact, and concentration compactness argument was used by Flucher and M\"uller \cite{FlucherMuller-exs} to proved existence of minimizers for general autonomous nonlinearity $f$ under the penalty condition $c<c_\infty$, which is realized, in terms of the Lagrangean density
$$
F(t)\eqdef\int_0^tf(s)ds,
$$
as $F(t)>F_\infty(t)$, where $F_\infty$ is the appropriate asymptotic counterpart of $F$.

In the case of subcritical nonlinearity, more refined realizations of the penalty condition $c<c_\infty$ are given by Sirakov \cite{Sirakov}. 
The best known case when the condition $F>F_\infty$ does not yield $c<c_\infty$ is the Brezis--Nirenberg problem on bounded domains with $N=3$. The reason to it is that the relevant asymptotic problem is a problem in $\R^N$ rather than in the original domain and multiplication of the solution of the asymptotic problem with a cut--off function when $N>3$ gives an unredeemable error, indeed, \cite{BN} provides a non-existence counterexample.
\subsection{Problems without homogeneity}
In problems without homogeneity, Lagrange multipliers for different levels of the constraint functional have to be evaluated, which remains in general an open problem. One can still consider in this case minimization of the functional
\be
\label{Func0}
G(u)=\frac{1}{2}\int_{\Omega}(|\nabla u|^2+\lambda u^2)dx-\int_\Omega F(x,u(x))dx,
\ee
whose Euler--Lagrange equation is \eqref{Theeq}, under the Nehari constraint 
$(G'(u),u)=0$, but the derivative of the constraint can be eliminated from the Euler--Lagrange equation only under additional assumptions. A typical sufficient condition that allows to use the Nehari constraint is increasing $s\mapsto f(x,s)/s$. We refer the reader to a typical result on this lines presented in Chapter 4 in the book of Willem \cite{Willem}). Nehari constraint approach will be not considered here. 
 
Insofar as one can verify the Palais--Smale condition for bounded critical sequences, when the nonlinearity lacks homogeneity, one has to consider unconstrained minimax statements of the mountain pass type.

Below we use the following terminology. A number $c\in\R$ is called a {\em critical value} of a $C^1$- functional $G$ in a Banach space, if it is the value $G(u)$ at a critical point $u$ ($G'(u)=0$), and it is called a {\em critical level} of a functional if there is a sequence $u_k$, called critical sequence, such that $G(u_k)\to c$ and $G'(u_k)\to 0$. The Palais--Smale condition ${\mathrm(PS)_c}$ says that every critical sequence at the level $c$ has a convergent subsequence, in which case the critical level becomes a critical value. 
 
There is a general result on weak convergence of critical sequences for semilinear elliptic equations to a non-zero solution due to Rabinowitz \cite{Rabinowitz}, but his proof does not verify ${\mathrm (PS)_c}$ condition.

There are several important reasons to verify ${\mathrm (PS)_c}$, such as: utility of minimax statements to calculate or estimate critical values; sorting obtained by different methods by the values of the functional; estimating the Morse index of a critical point; and existence of multiple critical points via minimax principles that require the ${\mathrm (PS)_c}$. 
On the other hand, the ${\mathrm (PS)_c}$-condition in non-compact problems is well known to fail on critical levels produced by divergent bounded sequences of the form 
\be
\label{bumps}
u_k=w+\sum_{n=1}^m g_k^{(n)}w_\infty,
\ee
where $g_k^{(n)}$ are pairwise asymptotically orthogonal sequences of transformations responsible for the loss of compactness (actions of translations or dilations), $w$ is a critical point of the functional and $w_\infty$ is a critical point of the asymptotic problem. 

Our approach is based on the observation that \eqref{bumps} is essentially the only way a bounded critical sequence can diverge, which leads to a conclusion that ${\mathrm (PS)_c}$ in the Ambrosetti-Rabinowitz settings holds whenever the mountain pass level satisfies $c<c_\infty$ (while $c\le c_\infty$ is true in general).
Decompositions of critical sequences similar to \eqref{bumps} have been introduced by Struwe \cite{Struwe} for the critical exponent case in bounded domains, followed by Brezis and Coron \cite{BrezisCoron} and by Lions \cite{Lions86} for subcritical problems in $\R^N$; see also Cao and Peng \cite{Cao} for the case of critical exponent on $\R^N$. We use a somewhat more detailed version of the ``multibump" decomposition, Theorem~\ref{abstractcc*} (\cite[Theorem~5.1]{ccbook}), which is a particular case of the abstract weak convergence decomposition from \cite{acc}. 

Verification of Palais--Smale condition at the mountain pass level is trivial when the problem has homogeneity and thus has an equivalent constrained minimization statement. This has been already observed by Cerami, Fortunato and Struwe \cite{CFS} who considered autonomous problem with the critical stem nonlinearity on bounded domains in the mountain pass setting, noting that the solvability condition $c<c_\infty$, established by Lions, extends to the mountain pass problems. Chabrowski and Yang \cite{ChabYang} have verified the Palais--Smale condition in the subcritical case under  an additional assumption $f_\infty(s)/s$ increasing for an interval of critical values $(0,J^\infty)$, where $J^\infty$ is a constrained minimum of $G_\infty$ under the Nehari constraint $(G'_\infty(u),u)=0$. Lemma~2.2 in \cite{PinoFelmer} states that $c_\infty$ is attained on the straight line path $t\mapsto tw$, where $w$ is, in terms of Chabrowski and Yang, a minimizer for $G_\infty$ under the Nehari constraint, which implies that $J^\infty=c_\infty$. In other words, the solvability condition of Chabrowski and Yang is an implicit form of the condition $c<c_\infty$, used in this paper. We extend (in Section 5) their result to the case when monotonicity of $f_\infty(s)/s$, assumed by Chabrowski and Yang, is no longer required.
Another related result for the subcritical case is due to Bartsch and Wang \cite{BartschWang}, but it is out of scope of this paper, as it deals with the equation $-\Delta u+V(x)u=f(x,u)$ for general subcritical $f$, deriving the Palais-Smale condition from unbounded (not necessarily coercive) $V>0$, namely such that for every $M>0$, the set
$V^{-1}(0,M)$ has a finite measure. 
We do not survey the literature here for the critical case, referring the reader to the bibliography in the books of Chabrowski, \cite{Chabrowski}, Flucher \cite{Flucher-book} and Willem \cite{Willem} and to the recent survey of Bartsch, Wang and Willem \cite{BWW}.
As a rule, the nonlinearity considered in literature is of the form $a(x)|u|^{2^*}$ plus a subcritical term, and this paper considers a more general case.

\vskip 5mm
Our paper is organized as follows. 
For the sake of simplicity, we assume that the nonlinearity $F(x,s)$ is continuously differentiable with respect to $s$ and that its derivative admits required asymptotic functions as uniform limits, so the asymptotic functionals and the asymptotic equations are well defined. 
Section 2 gives a generalization of the existence result of Flucher and M\"uller in the sense that the nonlinearity at infinity is defined due to discrete dilations, which generally gives a smaller $F_\infty$ (involved in the penalty condition $F>F_\infty$) than the upper limit defined in \cite{FlucherMuller-exs}.  A case in point here is the nonlinearity such that $t\mapsto F(e^t)e^{-2^*t}$ is a periodic function. This periodicity implies that $F(s)$ oscillates about the critical ``stem'' $s^{2^*}$. Some sort of oscillations about the critical ``stem'' are necessary for existence of solutions in the zero mass case ($\lambda=0$). Specifically, under the growth bounds \eqref{hom_sub} in the zero mass case, the mountain pass solutions of the autonomous equation \eqref{Theeq} on $\R^N$ (which are equivalently provided as constrained minima)  satisfy the well-known Poho\v{z}aev identity (see \cite{Poh} for bounded domain, \cite{BerestickyLions} for $\R^N$): 
\be
\label{Poh}
\int_{\R^N} |\nabla u|^2dx=2^*\int_{\R^N} (F(u)-\lambda u^2)dx.
\ee
Validity of this identity requires that $u$ and its gradient are decaying sufficiently fast at infinity. These decay rates are verified in the case $\lambda>0$ with subcritical nonlinearity in \cite{BerestickyLions} and for $\lambda=0$ and the nonlinearity $F$ bounded by the critical stem $C|u|^{2^*}$ in \cite{FluCherMuller-decay}. 
For positive solutions, Poho\v{z}aev identity is equivalent to
$$
\int_0^\infty s^{2^*+1}\dfrac{d}{ds}\frac{F(s)-\lambda s^2}{s^{2^*}} d|{\{u\le s\}}|=0,
$$
which implies that $\frac{F(s)-\lambda s^2}{s^{2^*}}$ is necessarily non-monotone, unless it is a constant.  If $\lambda>0$, this relation is satisfied whenever $F(s)=o(s^2)$ at zero and $F(s)/s^2\to\infty$ at infinity, which are typical sufficient conditions for continuity of $G$ to have the mountain pass geometry. For $\lambda=0$, this condition becomes, however, a significant condition of oscillatory behavior of $F(s)/s^{2^*}$. 

Critical points obtained in Section 2 are used in calculations of Section 3 dealing with the non-autonomous critical problem in $\R^N$ in the zero mass ($\lambda=0$) case.  In Section 4 we deal with critical problems on bounded domains, that is, with a generalization of the Brezis--Nirenberg problem to oscillatory critical nonlinearities. In Section 5 we give an improved existence condition in the subcritical case on $\R^N$. 
The main results of this paper are existence theorems Theorem~\ref{main*}, Theorem~\ref{main**} and Theorem~\ref{main-sub}. For the sake of consistency we also include for each of these three cases an elementary statement that the non-strict version $c\le c_\infty$ of the penalty condition is unconditionally true.
\section{Autonomous critical problem}
We consider the space $\merhav$, $N>2$, a completion of $C_0^\infty(\Omega)$ with respect to the gradient norm
$$
\|u\|=\left(\int_{\R^N}|\nabla u|^2dx\right)^{1/2},
$$
and we equip the space $\merhav$ with the group $D(N,\Z,\gamma)$ of unitary operators generated by the shifts
$$
D_{\R^N}\eqdef \{u\mapsto u(\cdot-y), y\in\R^N\}
$$
and by the action of discrete dilations with a fixed scaling factor $\gamma>1$,
$$
\delta_{\Z,\gamma}\eqdef\{u\mapsto \gamma^{\frac{N-2}{2}j}u(\gamma^j\cdot), j\in\Z\}. 
$$
Let $f\in C(\R)$ and let
$$
F(s)\eqdef\int_0^sf(\sigma)d\sigma.
$$
Assuming 
\be
\label{prel_assumptions}
|F(s)|\le C|s|^{2^*},%\mbox{ and } \lim_{s\to+\infty}F(s)/s^2=+\infty
\ee 
we set 
\be
\label{psi}
\psi(u)\eqdef\int_{\R^N}F(u)dx,
\ee
and 
\be
\label{thefunctional}
G(u)=\eqdef\frac12\|u^2\|- \psi(u),
\ee
and note that $G\in C(\merhav)$.
Let 
\be
\label{thepaths}
\Phi_G=\{\varphi\in C([0,\infty)\to H^1(\R^N)): \varphi_0=0, \lim_{t\to\infty}G(\varphi_t)=-\infty\}.\ee
If $\Phi_G=\emptyset$, we set $c(G)=+\infty$, otherwise 
\be
\label{critvalue}
c(G)\eqdef\inf_{\varphi\in\Phi_G}\sup_{t\in[0,\infty)} G(\varphi_t).
\ee.
\begin{proposition}
\label{emptyset} The set $\Phi_G$ is nonempty if and only if
$\sup F>0$, in which case $\Phi_G$ contains a path 
$u_t(x)=u(x/t)$ with $u\in C_0^\infty$ and $\psi(u)>0$.  
\end{proposition}
\begin{proof}
Not first that the path $u_s$ is a continuous map from $(0,\infty)$ to $\merhav$ that extends by continuity as $u_0=0$: $\|u_s\|^2=s^{N-2}\|u\|^2$, and in particular $\lim_{s\to 0}\|u_s\|=0$. Since the norm is continuous, it suffices to prove continuity of $u_s$ in $\mathcal D'(\R^N)$. Indeed, with $\varphi\in\core$ and $s,s_0>0$, by Lebesgue convergence theorem,
$$\int u_s\varphi=s^N\int u\varphi(s\cdot)\to  s_0^N\int u\varphi(s_0\cdot)\mbox{ as } s\to s_0.$$ 

If $\sup F>0$, then $\sup\psi>0$ and $\Phi_G\neq \emptyset$, 
since it contains a path $u_s(x)=u(x/s)$ with $\psi(u)>0$. 
Indeed, by the change of variables in respective integrals,
$$
G(u_s)=\frac12s^{N-2}\|u\|^2-s^N\psi(u).
$$
If, on the other hand, $F\le 0$, then $G\ge 0$ and $\Psi_G=\emptyset$.
\end{proof}
We assume that the following limits exist:
\bean
f_+(s)&\eqdef&\lim_{j\in\Z, j\to+\infty} \gamma^{\frac{N+2}{2}j}f(\gamma^{-\frac{N-2}{2}j}s),\\
f_-(s)&\eqdef&\lim_{j\in\Z, j\to-\infty} \gamma^{\frac{N+2}{2}j}f(\gamma^{-\frac{N-2}{2}j}s)
\eean
Repeating the definitions above for the functions $f_\pm$, we consider $F_\pm$,
$\psi_\pm$ $G_\pm$, $\Phi_{G_{\pm}}$ and $c(G_{\pm})$.
We will call the function $F$ {\em selfsimilar} with a factor $\gamma$ if
\be
\label{fpm}
F(s)=\gamma^{-Nj}F(s)(\gamma^{\frac{N-2}{2}j}s), j\in\Z, s\in\R.
\ee
It is uniquely defined by its values on the intervals $(1,\gamma)$
and $(-\gamma,-1)$. 
If $F$ is differentiable, one obviously has
\be
f(s)=\gamma^{-\frac{N+2}{2}j}f(s)(\gamma^{\frac{N-2}{2}j}s), j\in\Z, s\in\R.
\ee
Note that for any given $F$ that admits asymptotic functions $F_\pm$, they are selfsimilar.

The following statement generalizes Theorem~5.2 from \cite{ccbook}.
\begin{proposition}
\label{constraint}
Assume \eqref{prel_assumptions} and assume , for each of the signs ``$+$" and ``$-$",
that either $\sup F_\pm>0$, or $F_\pm=0$ with $\sup F>0$. Let
\be
\label{kappa}
\kappa(t)\eqdef\sup_{\|u\|^2=t}\psi(u), \,t>0
\ee
and 
Let $\kappa_\pm(t)$ be the value \eqref{kappa} corresponding to the functionals $\psi_\pm$.
If $F$ satisfies \eqref{fpm}, or if $\kappa(1)>\max\{\kappa_-(1),\kappa_+(1)\}$, then 
the maximum in \eqref{kappa} is attained.

Furthemore, the inequality $\kappa(t)>\kappa_\pm(t)$ holds whenever $F\ge F_\pm$ with the strict inequality in a neighborhound of zero.
\end{proposition}
\begin{proof} 

1. By substitution $u(s)=v(s/t^\frac{1}{N-2})$
\be
\label{kappa-t}
\kappa(t)=\kappa(1)t^\frac{2^*}{2},
\ee
so it suffices to prove the lemma for $t=1$. 

2. Assume now that $F$ satisfies \eqref{fpm}.
If $F=0$, then $\kappa=0$ and any function with the given norm is a maximizer for $\psi$. Assume now that $\sup F>0$. 
Let $u_k$ be a minimizing sequence in \eqref{kappa}, that is,
$\|u_k\|^2=1$ and $\psi(u_k)\to\kappa(1)$. Let $y_k^{(n)}\in\R^N$, $j_k^{(n)}\in\Z$, $w^{(n)}\in\merhav$ and the index sets $\N_{+\infty},\N_{-\infty},\N_{0}\subset\N$ be as in Theorem~\ref{abstractcc*}. Note that by Lemma~\ref{BL}, 
\be
\label{mass0}
\kappa(t)=\lim\psi(u_k)=\sum_{n\in\N}\psi(w^{(n)}).
\ee
At the same time from \eqref{norms*}
\be
\label{nrg0}
1=\|u_k\|^2\ge\sum_{n\in\N}\|w^{(n)}\|^2.
\ee
Let $v^{(n)}(x)=w^{(n)}(s_nx)$ with 
$s_n=\|w^{(n)}\|^\frac{2}{N-2}$, which gives whenever $w^{(n)}\neq 0$, $\|v^{(n)}\|^2=1$, otherwise $v^{(n)}=0$.
Then \eqref{mass0}, restated in the terms of $v^{(n)}$, holds
\be
\kappa(1)=\sum_{n\in\N}s_n^N\psi(v^{(n)})\le \kappa(1)\sum_{n\in\N}s_n^N,
\ee
or, in other words,
\bea
\label{mass01}
\sum_{n\in\N}s_n^N\ge 1.
\eea
On the other hand, from \eqref{nrg0} we derive
\be
\label{nrg01}
\sum_{n\in\N}s_n^{N-2}\le 1.
\ee
Relations \eqref{mass01} and \eqref{nrg01} can hold simultaneously if and only if there is a $n_0\in\N$
such that $s_{n_0}=1$, while $s_n=0$ whenever $n\neq n_0$.
Consequently, we have from \eqref{BBasymptotics*}
$u_k-\gamma^{\frac{N-2}{2}j_k^{(n_0)}} w^{(n_0)}(\gamma^{j_k^{(n_0)}}(\cdot-y_k^{(n_0)}))\to 0$ in $L^{2^*}(\R^N)$,
or, by replacing the maximizing sequence $u_k$ with the maximizing sequence
$\hat u_k\eqdef \gamma^{-\frac{N-2}{2}j_k^{(n_0)}} u_k(\gamma^{-j_k^{(n_0)}}\cdot+y_k^{(n_0)})$, we have
$\hat u_k\to w^{(n_0)}$ in $L^{2^*}(\R^N)$. At the same time,
by the weak lower semicontinuity of norms, $\|w^{(n_0)}\|^2\le 1$, while
by continuity of $\psi$ in $L^{2^*}(\R^N)$ we have $\psi(w^{(n_0)})=\kappa(1)$.
Since the latter can hold only when $\|w^{(n_0)}\|^2= 1$, $w^{(n_0)}$ is the desired maximizer.

3. Consider now the general case with $\kappa(1)>\kappa_\pm(1)$ and note that, since $F_\pm$
satisfies \eqref{fpm}, the maximum for $\kappa_\pm(1)$ is attained due to the previous step. Let $u_k$ be a minimizing sequence in \eqref{kappa} and let $y_k^{n}\in\R^N$, $j_k^{n}\in\Z$, $w^{n}\in\merhav$ and the index sets $\N_{+\infty},\N_{-\infty},\N_{0}\subset\N$ be as in Theorem~\ref{abstractcc*}. Note that by Lemma~\ref{BL},
\be
\label{mass02}
\kappa(t)=\lim\psi(u_k)=\sum_{n\in\N_0}\psi(w^{(n)})+\sum_{n\in\N_{-\infty}}\psi_-(w^{(n)})+\sum_{n\in\N_{+\infty}}\psi_+(w^{(n)}),
\ee
Let, as in the step 2, $v^{(n)}(x)=w^{(n)}(s_nx)$ with $s_n=\|w^{(n)}\|^\frac{2}{N-2}$, which yields 
\bea
\label{mass03}
\sum_{n\in\N_0}s_n^N & + & \frac{\kappa_-(1)}{\kappa(1)}\sum_{n\in\N_{-\infty}}s_n^N\\
&+&\frac{\kappa_+(1)}{\kappa(1)}\sum_{n\in\N_{+\infty}}s_n^N\ge 1.
\eea
while we still have \eqref{nrg01}.
Since $\frac{\kappa_\pm(1}{\kappa(1)})<1$, relations \eqref{mass03} and \eqref{nrg01} can hold simultaneously if and only if there is a $n_0\in\N_0$
such that $s_{n_0}=1$, while $s_n=0$ whenever $n\neq n_0$.
Consequently, we have from \eqref{BBasymptotics*}
$u_k-w^{(n_0)}(\cdot-y_k^{(n_0)})\to 0$ in $L^{2^*}(\R^N)$.
Similarly to the step 2 we conclude that 
$w^{(n_0)}$ is the desired maximizer.

4. The last assertion of the proposition is obvious if we take into account that the range of any function in $\merhav$ is a connected set whose closure contains zero.
\end{proof}
\begin{remark}
\label{justrem}
Note that one always has $\kappa(t)\ge\kappa_\pm(t)$, since if $w$ is a maximizer for 
$\kappa_\pm(1)$, then $\psi(\gamma{\frac{N-2}{2}j}w(\gamma^{j}x))\to \psi_\pm(w)$ as $j\in\Z$, $j\to\pm\infty$.

If $w_1$ is a maximizer of \eqref{kappa} for $t=1$, then, obviously,the function
$$
w_t(x)=w_1(t^{-\frac{1}{N-2}}x)
$$
is a maximizer for $t>0$.  
%It is shown in \cite{FluCherMuller-decay} that the maximizer $w_t$ is, up to shifts, a %decreasing radial function, and under general conditions (\cite{CEF,PS1,PS2}, which includes %in particular $F(s)=|s|^{2^*}$) it is unique. 
\end{remark}
We now connect the maximizers \eqref{kappa} with the mountain pass values for $G$.
\begin{proposition} Assume that $\sup F>0$ and that
\label{exist-crit-infty} 
\be
\label{hom_sub}
|f(s)|\le C|s|^{2^*-1}.%\mbox{ and } \lim_{s\to+\infty}F(s)/s^2=+\infty
\ee 
Then
any maximizer $w$ for \eqref{kappa} corresponding to 
$$t=t_0\eqdef (2^*\kappa(1))^{-\frac{2}{N-2}},$$
is a critical point of $G$. Moreover, 
\begin{itemize}
\item[(i)] $0<G(w)=c(G)=\max_{t\ge 0}G(w(\cdot/t))$; and
\item[(ii)] if $v$ is a critical point of $G$ such that $\psi(v)>0$,
then $G(v)\ge c(G)$.
\end{itemize}
\end{proposition}
\begin{proof} From \eqref{hom_sub} follows that $G\in C^1(\merhav)$.
We prove (ii) first. Let $v$ be any critical point with $\psi(v)>0$ and let $v_s(x)=v(x/s)$. Since $\psi(v)>0$, we have $v_s\in\Phi_G$ and the function $s\mapsto G(v_s)=\frac12s^{N-2}\|v\|^2-s^N\psi(v)$  has a single critical point, a 
maximum, which is necessarily attained at $s=1$ since $G'(v)=0$. Then
$$c(G)\le \max_{s\ge 0} G(v_s)=G(v),$$
which verifies (ii).

Assume now that \eqref{kappa} has a maximizer $w$.
By Remark~\ref{justrem}, $w$ is a critical point of $G$.
By the argument above, $c(G)\le G(w)$.
On the other hand, since any path $u_s\in\Phi_G$ starts at the origin and is unbounded, we have, using the notation $\hat u=u(\|u\|^{-\frac{2}{N-2}}\cdot)$, $r_s=\|u_s\|$, and taking into account that $\|\hat u\|=1$,
\bean
c(G) &= &\inf_{u_s\in\Phi}\max_{s\ge 0}\frac12 r_s^2 - r_s^\frac{2N}{N-2}\psi(\hat u_s)\\
&\ge &\inf_{u_s\in\Phi}\max_{s\ge 0}\frac12 r_s^2 - r_s^\frac{2N}{N-2}\kappa(1)\\
&=&\max_{r\ge0}\frac12 r^2 - r^\frac{2N}{N-2}\kappa(1)=G(w).
\eean
We conclude that $c(G)=G(w)$
and the path $w_t\in \Phi_G$ is a minimal path passing through the critical point $w$.
Explicit calculation of maximum on the minimal path gives the value of $t_0$. This proves (i).
\end{proof}
%The following property is a minimax counterpart of a classical property of %constrained minima with lack of compactness, an inequality between the critical value %and the critical value of the asymptotic functional:
%\begin{proposition}    
%\label{c--c0} If $F_\pm$ is definied, then
%\be\label{c-c0}
%c(G)\le c(G_\pm).\ee
%\end{proposition}
%\begin{proof} 
%Without loss of generality consider the ``$-$" case.
%If $\Phi_{G_{-}}=\emptyset$, the proposition is tautological.
%Assume that $\Phi_{G_{-}}\neq\emptyset$. Then due to Proposition~\ref{emptyset} 
%and Proposition~\ref{exist-crit-infty}, the functional $G_-$ has a critical point $w$
%and the path $w(\cdot/s)$, $s\ge 0$, is a minimal path.
%Let $$
%w_j(x)\eqdef \gamma^{\frac{N-2}{2}j}u(\gamma^{j}x), j\in\N.
%$$
%Then for all $j>0$ sufficiently large, 
%$$\psi(w_j)=\int_{\R^N} \gamma^{-Nj}F(\gamma^{\frac{N-2}{2}j}w(x))dx>0,$$
%and thus $\Phi(G)\neq\emptyset$. 
%Then
%\bean
%c(G)&\le&\max_{t\ge 0}  G(w_j(\cdot/t))\to \max_{t\ge 0} G_{-}(w_j(\cdot/t))\\
%&=&\max_{t\ge 0} G_{-}(w(\cdot/t))=c(G_).
%\eean
%\end{proof}
%%%%%%%%%%%%%%%%%%%%%%%%%%%%%%%%%%%%%%%%%%%%%%%%%%%%%%%%%%%%%%%%%%%%%%%%%%%%%%
\section{Non-autonomous critical problem in $\R^N$.}
Let now $f(x,s)\in C(\R^N\times\R)$, $N>2$.

Assume that for some $\gamma>1$ the following limits exist and that the convergence is uniform:
$$
f_{0}(s)\eqdef\lim_{|x|\to \infty} f(x,s),
$$
$$
f_{-}(s)\eqdef\lim_{j\in\Z, j\to -\infty} \gamma^{\frac{N+2}{2}j}f(\gamma^{j}x,\gamma^{-\frac{N-2}{2}j}s),
$$
$$
f_{+}(s)\eqdef\lim_{j\in\Z, j\to +\infty} \gamma^{\frac{N+2}{2}j} f(\gamma^{j}x,\gamma^{-\frac{N-2}{2}j}s).
$$
Let
$$F(x,s)=\int_0^sf(x,\sigma)d\sigma,$$
and assume the inequality
\be
\label{sbcrtcl}
|f(x,s)|\le
C|s|^{2^*-1}, s\in\R, x\in\R^N.
\ee
Let
$$\psi(u)=\int_{\R^N}F(x,u)dx,$$
\be \label{the_func} G(u)=\frac12\|u\|^2_{\merhav}-\psi(u),\ee
and let $\Phi_G$, and $c(G)$ be as in \eqref{thepaths} and
in \eqref{critvalue} respectively.
We will also consider similarly defined respective $F_{\#}$, 
$\psi_{\#}$, $G_{\#}$, $\Phi_{G_\#}$ 
$\Phi_{G_\#}$ and $c(G_\#)$, 
where $\#$ will refer in what follows to any of the three indices $0$, $+$, $-$.
\begin{proposition} 
\label{c--c} Assume \eqref{sbcrtcl}.
%Let $\sup F_\#>0$ or ($F_\#=0$ and $\sup F>0$).
Then 
\be\label{c-c}
c(G)\le c(G_\#).\ee
\end{proposition}
\begin{proof}
If $F_\#\le 0$, then $\Phi_{G_\#}=\emptyset$,
$c(G_\#)=+\infty$ and the statement is tautological. We may assume now that
$\sup F_\#>0$. By Proposition~\ref{exist-crit-infty}, each of the functionals $G_\pm$ has a critical point $w_\pm$ on a minimal path $w_\pm(\cdot/s)$. If the functional $G_0$ does not have a critical point, from Proposition~\ref{constraint} and Proposition~\ref{exist-crit-infty} it follows that its associated constrained supremum $\kappa(t)$ equals to one of its asymptotic functionals $\psi_{\pm}$, which are in fact are given by the original $F_\pm$ and consequently, we have $c(G_0)0=c(G_+)$ or $c(G_0)=c(G_-)$. In this case the inequality  $c(G)<c(G_0)$ tautologically follows from the inequalities with $G_\#=G_\pm$ for which critical points do exist. Thus, without loss of generality we may assume that $G_\#$ has a critical point $w$ with a minimal path
$w(\cdot/s)\in\Phi_{G_\#}$.
Then there is a sequence of $j_k\in\Z$ and $y_k\in \R^N$
such that either $j_k\to +\infty$, or $j_k\to -\infty$, or $|y_k|\to\infty$, 
and such that, with
$$
w_{t;k}\eqdef\gamma^{\frac{N-2}{2}j_k}w_t(\gamma^{j_k}\cdot+y_k),
$$
the sequence $G(w_{t;k})$ converges to $G_\#(w_t)$ uniformly in $t$.
Indeed, with 
$$
F_k(x,s)=\gamma^{Nj_k}F(\gamma^{-j_k}(x-y_k),\gamma^{\frac{N-2}{2}j_k}s),
$$
it is easy to see that for any $\epsilon>0$ there is a $k_\epsilon\in\N$
so that for all $k>k_\epsilon$ and $t\ge \epsilon,$
$$
\psi(w_{t;k})=\int_{\R^N} F(x,w_{t;k})=t^{N}\int_{\R^N} F_k(tx,w)
\to t^{N}\psi_\#(w)=\psi_\#(w_t).
$$
Let us redefine the path $w_{t;k}$ for $t\in[0,\epsilon]$,
by $\frac{t}{\epsilon}w_{\epsilon;k}$. Then it is easy to see that with $\epsilon$ sufficiently small,  
\bean
\label{comparison}
c(G)&\le&\max_{t\ge \epsilon}  G(w_{t;k})\to \max_{t\ge 0} G_\#(w_t)=c(G_\#).
\eean
\end{proof}
\begin{theorem}
\label{main*}
% Let $\sup F_\#>0$ or $F_\#=0$ with $\sup F>0$. 
Assume in addition to \eqref{sbcrtcl} that $F$ is not identically $0$ and 
\begin{itemize}
%\item[(M)] the maximum in \eqref{kappa} for $\psi_0$ is attained -- which is true if $F_0$ satisfies conditions of Propositon~\ref{constraint} (in particular, if 
%$F_0\ge F\pm$ with a strict inequality  in a neighborhood of zero);
%and
\item[{\rm (R)}] 
there exists $\mu>2$, such that
$$
f(x,s)s\ge\mu F(x,s), s\in\R\setminus\{0\}, x\in\R^N.
$$
\end{itemize}
If 
\be\label{crin}c(G)<c(G_\#),\ee
then the functional\eqref{the_func} has a critical point at the level $c(G)$ and every 
critical sequence at this level has a convergent subsequence.
Furthermore, the relation \eqref{crin} is satisfied if
\be
\label{>}
F(x,s)\ge F_{\#}(s),  x\in\R^N,
\mbox{with the strict inequality for $s$ in a neighborhood of zero.}
\ee
\end{theorem}
\begin{proof} Since $c(G)<\infty$, $\sup F>0$ and from the condition (R) trivially follows the mountain pass geometry, namely, the non-empty $\Psi_G$ and $c(G)>0$. Note also that one can pass to the limit in (R), so the condition holds also for $F_\#$.
By the standard mountain pass reasoning, the functional $G$ possesses a critical sequence $u_k\in\merhav$, that is, $G(u_k)\to c(G)$ and $G'(u_k)\to 0$. Also by a standard argument, it follows from (R) that the sequence $u_k$ is bounded in $\merhav$. Consider now the renamed subsequence of $u_k$ given by Theorem~\ref{abstractcc*}.
If $w^{(n)}= 0$ in \eqref{BBasymptotics*} for all $n\ge 2$, then 
 $u_k\to w^{(1)}$ in $L^{2^*}$, $\psi'(u_k)$ converges in $\merhav$, and from $G'(u_k)\to 0$ it follows that $u_k$  converges in $\merhav$ to a critical point of $G$ at the level $c(G)$.

Let us assume now that for some $m\ge 2$, $w^{(m)}\neq 0$.
Due to \eqref{norms*} and Lemma~\ref{ser**pen}, we have the following estimate of $c(G)$ from below:
\be
\label{frombelow}
c(G)=\lim G(u_k)\ge G(w^{(1)})+\sum_{n\in\N_{0}}G_0(w^{(n)})+\sum_{n\in\N_{+\infty}}G_{+}(w^{(n)})+\sum_{n\in\N_{-\infty}}G_{-}(w^{(n)}).
\ee
Note that with necessity, $w^{(1)}$ is a critical point of $G$, and $w^{(n)}$, $n\ge 2$, are critical points of correspondent $G_{\#}$. 
Let $G_m=G_0$ if $m\in\N_0$, $G_m=G_+$ if $m\in\N_+$ and $G_m=G_-$ if $m\in\N_-$. The correspondent asymptotic nonlinearity we will denote as $F_m$. Due to condition (R),
$G(w^{(1)})=\int [\frac12 f(x,w^{(1)})w^{(1)}-F(x,w^{(1)})]\ge 0$.
Similarly $G_\#(w^{(n)})\ge 0$. Furthermore, since $w^{(m)}\neq 0$, $\sup F_m>0$
and then, for $n=m$, (R) implies $\psi_m(w^{(m)})>0$ and $G_m(w^{(m)})> 0$.    
Combining \eqref{frombelow} with \eqref{crin}, we have
$$
G_\#(w^{(m)})<c(G_\#).
$$
On the other hand, by Proposition~\ref{exist-crit-infty}, $G_\#(w^{(m)})\ge c(G_\#)$,
which is a contradiction.
Thus our assumption above that there is such $m\ge 2$ is false. Consequently,  $w^{(1)}$ is the desired critical point of $G$.
\par
It remains to verify that \eqref{>} implies \eqref{crin}. If $G_{\#}$ has a critical point $w$ with the critical value $c(G_\#)$ lying on a minimal path $w_t(x)=w(x/s)$ (due to Proposition~\ref{constraint} this is always the case when $\#$ is a $+$ or a $-$), then from \eqref{>} it follows that $\psi(w)>\psi_\#(w)$, and then
$c(G)\le \max_{s\ge 0} G(w_s)<\max_{s\ge 0} G_\#(w_s)=c(G_\#)$.
Assume now that the maximum in \eqref{kappa} for $\psi_\#$ is not attained, which is, with necessity the case of $\psi_0$. Then, by Proposition~\ref{constraint} and Remark~\ref{justrem}, either
$\kappa_0(t)=\kappa_+(t)$ or $\kappa_0(t)=\kappa_-(t)$, which, by Proposition~\ref{exist-crit-infty} gives immediately $c(G_0)=c(G_+)$ or $c(G_0)=c(G_-)$. In either case the inequality $c(G)<c(G_0)$ follows from 
one of the two verified conditions $c(G)<c(G_+)$ or $c(G)<c(G_-)$.
$f$\end{proof}
\section{Critical case, problems in domains}
Let now $\Omega\subset\R^N$, $N>2$, be a domain with $\overline\Omega\neq \R^N$.
Let $\mathcal D^{1,2}_0(\Omega)$ be the closure of $C_0^\infty(\Omega)$ in $\merhav$. 
Note that $\mathcal D^{1,2}_0(\Omega)$ coincides with  $H^1_0(\Omega))$ when
$|\Omega|<\infty$, but generally the elements of this space are bounded only in a weighted $L^2$-norm. 

We consider here the functional \eqref{the_func}, restricted to $\mathcal D^{1,2}_0(\Omega)$.
We assume that $f(x,s)\in C(\R^N\times\R)$ satisfies  \ref{sbcrtcl} if $|\Omega|0\infty$ and that
\be
\label{sbcrtcl-bdd}
|f(x,s)|\le
C(1+|s|^{2^*-1}), s\in\R, x\in\R^N
\ee
if $|\Omega|<\infty$. 
We consider the mountain pass problem 
$$
c(G;\Omega)\eqdef\inf_{\varphi\in\Phi_{G,\Omega}}\sup_{t\ge 0}G(\varphi_t),
$$
where $\Phi_{G,\Omega}\subset \Phi_{G,\Omega}$ consists of paths with values in $\mathcal D^{1,2}_0(\Omega)$.
\begin{remark}
\label{rem:bdd}
Regarding $\mathcal D^{1,2}_0(\Omega)$ as a closed subspace of $\merhav$,
one can apply Theorem~\ref{abstractcc*} to bounded sequences $u_k\in \mathcal D^{1,2}_0(\Omega)\subset \merhav$. The weak limits $w^{(n)}$ 
 are of course not necessarily supported in $\Omega$. 
If, moreover, $|\Omega|<\infty$, then 
\begin{itemize}
\item[(a)] The set $\N_0$ consists only of the index $1$, since $u_k(\cdot-y_k)\rightharpoonup 0$ whenever $|y_k|\to\infty$. 
\item[(b)] The set $\N_{-\infty}$ is empty since otherwise it easily follows that $\|u_k\|_{L^2}\to\infty$.
\end{itemize}
\end{remark}
%Note also that the limit function $f_+$, being selfsimilar, satisfies \eqref{sbcrtcl}, so %that $G_+\in C^1(\merhav)$.
\begin{theorem}
\label{main**} Let $|\Omega|<\infty$ Assume in addition to \eqref{sbcrtcl-bdd} 
the condition (R) (restricted to $x\in\Omega$).
If 
\be\label{crin-bdd}c(G;\Omega)<c(G_+),\ee
then the functional\eqref{the_func} has a critical point at the level $c(G;\Omega)$ and every 
critical sequence at this level has a convergent subsequence.
\end{theorem}
\begin{proof} The argument is completely analogous to the proof of Theorem~\ref{main*}, with the critical sequence regarded as a sequence in $\merhav$, and with obvious simplifications due to Remark~\ref{rem:bdd} which leaves $+$ as the only value for $\#$.
\end{proof}
\begin{remark}
\label{main--}
Repetition of the proof of Theorem~\ref{main*} also gives that if $F$ is like in Theorem~\ref{main*}, $\Omega\subset\R^N$ is a domain and inequalities 
\eqref{crin} are satisfied, then the functional $G$ understood in restriction to $\mathcal D^{1,2}_0(\Omega)$ satisfies $(PS)_c$-condition at the mountain pass level $c(G;\Omega)$
and thus has a critical point at this level. 
\end{remark}

\begin{remark}
Theorem~\ref{main**} and Remark~\ref{main--}, unlike their counterpart on $\R^N$, cannot claim that $c(G;\Omega)<c(G_\#)$ follows from $F(x,s)>F_\#(s)$, since the problem at infinity is supported on a different domain. Furthermore, \cite{BN} offers a counterexample
for $N=3$ and $F(s)=|s|^{2^*}+\lambda u^2$ with sufficiently small positive $\lambda$. 

On the other hand, if $\Omega$ is a bounded domain, $N>3$ and $F(x,s)\ge F_+(s)+\epsilon s^2$ with some $\epsilon>0$, then 
$c(G;\Omega)<c(G_+)$. This follows by repetition of the estimates of Brezis and Nirenberg from \cite{BN} that involves the estimates for the minimizer of $\kappa_+$ in \cite{FluCherMuller-decay} by the Talenti minimizer.  
\end{remark}
\begin{proposition}
Under conditions of Theorem~\ref{main**},
$c(G;\Omega)\le c(G_+)$.
\end{proposition}
\begin{proof}
Let $w_k\in\core$ be a sequence convergent in $\merhav$ to a critical point $w$ for $G_+$ satisfying $G_+(w)=c(G_+)$.  
$$
w_{s;j,k}(x)=\gamma^{\frac{N-2}{2}j}w_k(\gamma^{j}x/s), j\in\N. 
$$
Assume without loss of generality that $0\in\Omega$. Than for any $\epsilon>0$, $s\ge\epsilon$, and all $j$ sufficiently large, $w_{s;j,k}\in\core$.
For $s<\epsilon$ we can redefine $w_{s;j,k}$ as $\frac{s}{\epsilon}w_{\epsilon;j,k}$.
Then, with 
$$
F_t(x,u)=t^{N}F(tx,t^{-\frac{N-2}{2}}u)\to F_+(u) \mbox{ as } t\to 0.
$$
and assuming, for each $k$, $\epsilon=\epsilon(k)$ sufficiently small,
we have, similarly to the argument of Proposition~\ref{c--c},
\bean
c(G)\le \max_{s\ge 0} G(w_{s;j,k})\to
\max_{s\ge 0} \frac12 s^{N-2}\|w_k\|^2-s^N\psi_+(w_k).
\eean
It is easy to see that as $k\to\infty$, the right hand side converges to 
$$\max_{s\ge 0} \frac12 s^{N-2}\|w\|^2-s^N\psi_+(w)=c(G_+).$$
\end{proof}
It is also easy to see that $c(G;\Omega)\le c(G_\#)$ for general $\Omega$.
\section{Subcritical case}
Consider now the ``positive mass'' case $\lambda>0$ with the functional $\psi$ defined by expression \eqref{psi} 
on the Sobolev space $H^1(\R^N)$, $N\ge 1$, equipped with the equivalent norm 
$$\|u\|^2=\int_{\R^N}(\nabla|u|+\lambda|u|^2)dx.$$
Let
\be
\label{SubFunc}
G(u)=\frac12\|u\|^2_{H^1(\R^N)}-\psi(u).
\ee
We assume for every $\epsilon>0$ there exist $p_\epsilon\in(2,2^*)$ and
$C_\epsilon>0$ such that
\be
\label{subcr-subcr}
|f(x,s)|\le
\epsilon(|s|+|s|^{2^*-1})+C_\epsilon|s|^{p_\epsilon-1}, s\in\R,
x\in\R^N,
\ee
or, for $N=1,2$,
\be
\label{subcr-subcr-12}
|f(x,s)|\le
\epsilon|s|+C_\epsilon|s|^{p_\epsilon-1}, s\in\R,
x\in\R^N,
\ee
which assures that $G\in C^1(H^1(\R^N))$.
Assuming that the following uniform limit exists,
$$f_\infty(s)\eqdef \lim_{|x|\to\infty}f(x,s),$$
we define by analogy $F_\infty$, $\psi_\infty$, $G_\infty$, etc.
\begin{proposition}
Assume that \eqref{subcr-subcr} holds.
Then
$$
c(G)\le c(G_\infty).
$$
\end{proposition}
\begin{proof}
Let $u_t\in\Phi_{G_+}$ and let $y_k\in\R^N$, $|y_k|\to\infty$.
Then 
$$
c(G)\le \max_{t\ge 0}G(u_t(\cdot-y_k))\to \max_{t\ge 0}G_\infty(u_t).
$$
Minimizing the inequalities over all paths in $\Phi_{G_+}$ we arrive at $c(G)\le c(G_\infty)$.
\end{proof}
\begin{theorem}
\label{main-sub} Assume \eqref{subcr-subcr} and (R).
If $c(G)<c(G_\infty)$, then every sequence $u_k\in H^1(\R^N)$,
such that $G'(u_k)\to 0$ and $G(u_k)\to c(G)$, has a subsequence convergent to a critical point of $G$.
Furthermore, the relation $c<c_\infty$ holds if 
\be
\label{**}
F(x,s)\ge F_\infty(s), x\in\R^N, \mbox{ with the strong inequality in a neighborhood of }s=0.
\ee

\end{theorem}
\begin{proof}  
The beginning of the proof is completely analogous to that for Theorem~\ref{main*} and can be abbreviated. We apply Theorem~\ref{abstractcc*} to the bounded critical sequence, noting, similarly to Remark~\ref{rem:bdd}, that $\N_{-\infty}=\emptyset$ and, moreover, $\N_{+\infty}=\emptyset$ since $F_+=0$. With \eqref{ser*pen} taking the role of \eqref{ser**pen} we arrive at an immediate analog of \eqref{frombelow}:
\be
\label{lemalra}
G(w^{(1)})+\sum_{n\ge 2}G_\infty(w^{(n)})\le c(G)<c(G_\infty).
\ee
As in the proof of Theorem~\ref{main*}, all the terms in the left hand side are non-negative due to (R). 
Assume that there exists $m\ge 2$ such that $w^{(m)}\neq 0$.
Observe the path $t\mapsto w^{(m)}(\cdot/t)$ is of the class $\Phi_{G_\infty}$.
Indeed, its continuity in $\merhav$-norm was shown in Proposition~\ref{emptyset}, and the proof of remaining continuity, in $L^2$ norm is analogous. The nonlinearity satisfies the requirements of \cite{BerestickyLions} for Poho\v{z}aev identity \eqref{Poh}. Thus, since
$$\int_{\R^N}(F_\infty(w^{(m)})-\lambda |w^{(m)}|^2)dx=\frac{1}{2^*}\|w^{(m)}\|_{\merhav}^2>0$$
the functional $G_\infty$ on the path $t\mapsto w^{(m)}(\cdot/t)$,
$$G_\infty(w^{(m)}(\cdot/t))= \frac12 t^{N-2}\|w^{(m)}\|_{\merhav}^2-t^N\int_{\R^N}(F_\infty(w^{(m)})-\lambda |w^{(m)}|^2)dx,$$
converges to $-\infty$ when $t\to\infty$.
Furthermore, the maximum of this expression over $t\ge 0$ is is clearly attained at a single point, which is necessarily $t=1$, since 
since $(G_\infty'(w^{(m)}),w^{(m)})=0$.
Thus 
\be
\label{cntr}
c(G_\infty)\le \max_{t\ge 0}G_\infty(w^{(m)}(\cdot/t))= 
G_\infty(w^{(m)})dx.
\ee 
On the other hand, since, $G_\infty(w^{(n)})\ge 0$ for all $n>1$ and $G(w^{(1)})\ge 0$,
we have from \eqref{lemalra},
$$G_\infty(w^{(m)})\le c(G)<c(G_\infty),$$
arriving to a contradiction with \eqref{cntr} unless $w^{(n)}=0$ for all $n\ge 2$. This, via the usual convergence argument, yields $u_k\to w^{(1)}$.
The implication \eqref{**} $\Rightarrow c(G)<c(G_\infty)$ is immediate once we take into account that the range of any function in $H^1_0(\Omega)$ is connected and contains zero.
\end{proof}
\begin{remark}
Theorem~\ref{main-sub} can be trivially generalized to the case of periodic coefficients. Let $V\in L^\infty(\R^N)$, $V>0$, be $\Z^N$-periodic, that is satisfy 
$$V(x+y)=V(x) \mbox{for all } x\in\R^N, y\in \Z^N$$ 
We may equip $H^1(\R^N)$ with an equivalent norm
\be
\|u\|^2=\int_{\R^N}(|\nabla u|^2+V(x)u^2)dx,
\ee
Assume that there is a function $f_\infty\in C(\R^{N}\times \R)$, $\Z^N$-periodic in the first argument, such that $f(x+y_k,s)\to f_\infty(x,s)$ for any sequence $y_k\in \Z^N$, $|y_k|\to\infty$. Theorem~\ref{main-sub} remains true also under these modifications. 
\end{remark}

\section{Appendix}
\subsection{Weak convergence decomposition in $\merhav$}
The following theorem is Theorem~ from \cite{ccbook}, with the dilation factor $2$
replaced by general $\gamma$.
\begin{theorem}
\label{abstractcc*}  Let $u_{k}\in \mathcal{D}^{1,2}(\R^N)$, $N>2$, be a
bounded sequence. Let $\gamma>1$. There exist $w^{(n)}  \in
\mathcal{D}^{1,2}(\R^N)$, $y_{k} ^{(n)} \in \R^N$, $j_{k} ^{(n)}
\in \Z$ with $k,n \in\N$, and disjoint sets
$\N_0,\N_{+\infty},\N_{-\infty}\subset\N$, such that, for a
renumbered subsequence of $u_k$,
\begin{eqnarray}
\label{w_n*}
&&
w^{(n)}=\wlim\,
\gamma^{-\frac{N-2}{2}j_k^{(n)}}u_k(\gamma^{-j_k^{(n)}}\cdot+y_k^{(n)}),\;
n\in\N,
\\
%\end{equation}
%\begin{equation}
\label{separates*}
&&
 |j_{k}^{(n)} -  j_{k}
^{(m)}|+|\gamma^{j_k^{(n)}}(y_{k} ^{(n)} - y_{k} ^{(m)})|\to \infty
\mbox{ for } n \neq m,
%\end{equation}
%\begin{equation}
\\
\label{norms*}
&&
\sum_{n\in \N} \|w
^{(n)}\|_{\mathcal{D}^{1,2}}^2 \le
\limsup\|u_k\|_{\mathcal{D}^{1,2}}^2,
%\end{equation}
%\begin{equation}
\\
\label{BBasymptotics*}
&&
u_{k}  -  \sum_{n\in\N}
\gamma^{\frac{N-2}{2}j_k^{(n)}} w^{(n)}(\gamma^{j_k^{(n)}}(\cdot-y_k^{(n)}))
  \to 0\; \mbox{ in } L^{2^*}(\R^N),
\end{eqnarray}
and the series above converges uniformly in $k$.
\par Furthermore,
$1\in\N_0$, $y_{k} ^{(1)} =0$; $j_k^{(n)}=0$ whenever $n\in\N_0 $;
$j_k^{(n)}\to -\infty$ (resp. $j_k^{(n)}\to +\infty$) whenever
$n\in\N_{-\infty}$ (resp. $n\in\N_{+\infty}$); and $y_k^{(n)}=0$
whenever $|\gamma^{j_k^{(n)}}y_k^{(n)}|$ is bounded.
\end{theorem}
The following statements are, respectively, Lemma~5.6 (an elementary modification) and Remark~3.4 from \cite{ccbook}.
\begin{lemma}
\label{BL} Let $F\in
C(R^N\times\R)$ satisfy $|F(x,s|\le c|s|^{2^*}$, let $\gamma>0$, $N>2$, and assume that the following limits exist and are uniform in $x\in\R^N$:
\begin{eqnarray*}
F_+(s)\eqdef\lim_{j\in\Z, j\to+\infty} \gamma^{-NJ}F(\gamma^{-j}x,\gamma^{\frac{N-2}{2}}j),\\
F_-(s)\eqdef\lim_{j\in\Z, j\to-\infty} \gamma^{-NJ}F(\gamma^{-j}x,\gamma^{\frac{N-2}{2}}j),\\
F_0(s)\eqdef\lim_{|x|\to\infty} F(x,s).
\end{eqnarray*}
Let $u_k\in\mathcal{D}^{1,2}(\R^N)$, $w^{(n)}$, $y_k^{(n)}\in\R^N$ and
let $j_k^{(n)}\in\Z$, $\N_0,\N_{+\infty},\N_{-\infty}\subset\N$ be as provided by Theorem~\ref{abstractcc*}.
Then
\begin{eqnarray}
\label{ser**pen} &&\lim_{k \to \infty} \int_{\R^N} F(u_k)
\\&=&
\nonumber
\sum_{n\in\N_0} \int_{\R^N} F_0(w^{(n)})+\sum_{n\in\N_{+\infty}} \int_{\R^N} F_+(w^{(n)})+\sum_{n\in\N_{-\infty}} \int_{\R^N} F_-(w^{(n)}).
\end{eqnarray}
\end{lemma}

\begin{lemma}
\label{BL0} Let $F\in
C(R^N\times\R)$ satisfy \ref{subcr-subcr}, \ref{subcr-subcr-12}, and assume that the following uniform limit exists:
\begin{eqnarray*}
F_\infty(s)\eqdef\lim_{|x|\to\infty} F(x,s).
\end{eqnarray*}
Let $u_k\in H^{1}(\R^N)$, $w^{(n)}$, $y_k^{(n)}\in\R^N$ be as provided by Theorem~\ref{abstractcc*} with $\N_\pm=\emptyset$ .
Then
\begin{equation}
\label{ser*pen} \lim_{k \to \infty} \int_{\R^N} F(u_k)=
\sum_{n\in\N} \int_{\R^N} F_\infty(w^{(n)}).
\end{equation}
\end{lemma}

\begin{center}
{\bf Acknowledgments}
\end{center}
The author thanks Moshe Marcus and Ian Schindler for their encouraging remarks.
This paper was written as Visiting Professor at University of Toulouse 1 and the author expresses his gratitude to J.Fleckinger and the rest of the faculty at Ceremath for their warm hospitality.

\end{document}